\documentclass{article}
\usepackage{spconf,graphicx}
\usepackage{amssymb}
\usepackage[cmex10]{amsmath}
\usepackage{pgfplots, filecontents}
\usepgfplotslibrary{patchplots}
\usepackage[]{algorithm2e}
\usepackage{enumerate}
\usepackage{color}
\usepackage{epstopdf}
\graphicspath{{./fig/}}
\usepackage{cite}
\newcommand{\ve}{\mathbf}
\newcommand{\m}{\mathbf}
\newcommand{\vea}[1]{\mathbf{\underline{#1}}}
\newcommand{\ma}[1]{\mathbf{\underline{#1}}}

\newcommand{\re}[1]{\Re\left\{{#1}\right\}} 
\newcommand{\im}[1]{\Im\left\{{#1}\right\}}

\title{Best Widely Linear Unbiased Estimator for Real Valued Parameter Vectors}
\name{O. Lang and M. Huemer}
\address{Johannes Kepler University Linz\\
	Institute of Signal Processing\\
	Altenbergerstra{\ss}e 69, 4040 Linz, Austria}
\begin{document}
%\ninept
%
\maketitle
\begin{abstract}
For classical estimation with an underlying linear model the best linear unbiased estimator (BLUE) is usually utilized for estimating the deterministic but unknown parameter vector. In the case of real valued parameter vectors but complex valued measurement matrices and noise vectors, the BLUE results in complex valued estimates, introducing a systematic error. In recent years widely linear estimators have been investigated for complex estimation. In this work a novel widely linear classical estimator is derived which incorporates the knowledge that the parameter vector is real valued. The proposed estimator is unbiased in the classical sense and it outperforms the BLUE and the best widely linear unbiased estimator (BWLUE) in terms of the variances of the vector estimator's elements.  
\end{abstract}
%

% -------------------------------------------------------------------
% Introduction
% -------------------------------------------------------------------
\section{Introduction}
\label{sec:intro}

In this work we consider the task of estimating a real valued parameter vector based on complex valued measurements. The measurements are assumed to be connected with the parameters via the linear model
\begin{equation}
\ve{y} = \m{H} \ve{x} + \ve{n}, \label{equ:WB_Real001}
\end{equation}
where $\ve{x}\in \mathbb{R}^{N_\ve{x} \times 1}$ is a real valued parameter vector, $\ve{y}\in \mathbb{C}^{N_\ve{y} \times 1}$ is a complex valued measurement vector, $\m{H}\in \mathbb{C}^{N_\ve{y} \times N_\ve{x}}$ is the complex valued measurement matrix, and $\ve{n}\in \mathbb{C}^{N_\ve{y} \times 1}$ is a complex valued zero mean random noise vector. 

In a Bayesian interpretation, the real valued parameter vector is improper and the application of widely linear estimators is obvious. For the definition of propriety we refer to \cite{Complex-Valued-Signal-Processing-The-Proper-Way-to-Deal-With-Impropriety} and Sec.~\ref{sec:AugmentedForm}. A widely linear Bayesian estimator suitable for the described problem is the widely linear minimum mean square error (WLMMSE) estimator. The WLMMSE estimator requires first and second order statistics about  $\ve{x}$. For zero mean $\ve{x}$, the WLMMSE estimator is of the form
\begin{equation}
\hat{\ve{x}} = \m{E}\ve{y} + \m{F}\ve{y}^*, \label{equ:WB_Real002}
\end{equation}
where $(\cdot)^*$ denotes the complex conjugate. Isolating the $i^\text{th}$ row of \eqref{equ:WB_Real002} yields
\begin{equation}
\hat{x}_i = \ve{e}_i^H\ve{y} + \ve{f}_i^H\ve{y}^*, \label{equ:WB_Real003}
\end{equation}
where $\ve{e}_i^H$ and $\ve{f}_i^H$ are the $i^\text{th}$ rows of $\m{E}$ and $\m{F}$, respectively, and $(\cdot)^H$ denotes the conjugate transposition. For real valued parameter vectors, the WLMMSE estimator fulfills the property
\begin{equation}
\ve{e}_i^H = \ve{f}_i^T, \label{equ:WB_Real004}
\end{equation}
with $(\cdot)^T$ denoting transposition. Consequently, in this case the WLMMSE estimator produces real valued estimates.

Other estimators such as the Bayesian LMMSE estimator \cite{Kay-Est.Theory}, the best linear unbiased estimator (BLUE) \cite{Kay-Est.Theory}, or the best widely linear unbiased estimator (BWLUE) \cite{Statistical-Signal-Processing-of-Complex-Valued-Data-The-Theory-of-Improper-and-Noncircular-Signals} do not result in real valued estimates $\hat{x}_i$ for real parameters $x_i$ in general, an exception is the case $\m{H}$,  $\ve{x}$, $\ve{n} \in \mathbb{R}$. In fact, the classical estimators BLUE and BWLUE do not allow the utilization of statistics about $\ve{x}$. These classical estimators are employed, e.g., when no statistics about $\ve{x}$ are available or when it is desired that the estimator is unbiased in the classical sense, i.e. it shall fulfill $E_{\ve{n}}[\hat{\ve{x}}]=\ve{x}$, where the index denotes the averaging probability density function (PDF). Since the classical estimators BLUE and BWLUE applied on the model in \eqref{equ:WB_Real001} do not result in real valued estimates, a systematic error is introduced. 

In this work, a novel classical estimator is derived which allows to incorporate the fact that $\ve{x}$ is real valued. This estimator is of widely linear form and unbiased in the classical sense. The derivation is executed by minimizing the variance of the vector estimator's elements subject to an unbiasedness constraint. This unbiasedness constraint incorporates the fact that $\ve{x}$ is real valued. Furthermore, in order to yield real valued estimates, $\ve{e}_i^H = \ve{f}_i^T$ is enforced during the derivation. The resulting widely linear estimator is termed BWLUE \emph{for real valued parameter vectors} and it in general outperforms the BLUE and the BWLUE in terms of the variances of the vector estimator's elements. A special case is investigated where the proposed estimator coincides with the real part of the BLUE. The estimator's performance is validated with a simulation example. 

The remainder of this paper is organized as follows: In Sec.~\ref{sec:AugmentedForm}, the basic concepts necessary for applying widely linear estimators are briefly recapitulated. Sec.~\ref{sec:Derivation} contains the derivation of the BWLUE for real parameter vectors, where we focus on proper noise statistics. A discussion about the proposed estimator is presented in Sec.~\ref{sec:Discussion} and a simulation example in which the proposed estimator is compared with several other estimators is provided in Sec.~\ref{sec:Simulation Results}.

%In addition, another classical estimator is derived for the case when the noise is improper. For that case it is assumed that the noise covariances are real valued in order to keep the formulas simple. The most general case with improper noise and complex noise covariance matrices is handled in an upcoming paper. The resulting estimators are termed BWLUE \emph{for real parameter vectors} for both, proper and improper noise. 

%In the following, the augmented form and improper statistics are used. For an excellent introduction to the augmented form, proper/improper statistics and widely linear processing we refer to \cite{Schreier-Buch} and \cite{Schreier-2011}.

%\emph{Notation:}
%\\ 
%Lower-case bold face variables ($\ve{a}$, $\ve{b}$,...) indicate vectors, and upper-case bold face variables ($\m{A}$, $\m{B}$,...) indicate matrices. We further use $\mathbb{R}$ and $\mathbb{C}$ to denote the set of real and complex numbers, respectively, $(\cdot)^T$ to denote transposition, $(\cdot)^H$ to denote conjugate transposition, $(\cdot)^*$ to denote the conjugate and $E[\cdot]$ to denote the expectation operator. $\re{\cdot}$ and $\im{\cdot}$ denote the real and imaginary part of a complex variable, respectively. 

\section{Preliminaries for Widely Linear Estimators} \label{sec:AugmentedForm}
In this section we recapitulate the preliminaries required to apply the widely linear estimators used in this work. This section is more or less a shortened version of the corresponding parts in \cite{Complex-Valued-Signal-Processing-The-Proper-Way-to-Deal-With-Impropriety}.

We start by constructing the \textit{complex augmented vector} $\underline{\ve{a}}$ of a vector $\ve{a}\in\mathbb{C}^{N_\ve{a}\times 1}$ by stacking $\ve{a}$ on top of its complex conjugate $\ve{a}^{*}$, i.e.
\begin{equation}
	\underline{\ve{a}} = \begin{bmatrix} \ve{a} \\ \ \ve{a}^{*} \end{bmatrix}\in\mathbb{C}^{2N_\ve{a}\times 1}.
\end{equation}
In order to characterize the second-order statistical properties of $\vea{a}$ we consider the augmented covariance matrix 
\begin{align}
	\underline{\m{C}}_{\ve{a}\ve{a}} 
			&= E[(\underline{\ve{a}} - E[\underline{\ve{a}}])(\underline{\ve{a}} - E[\underline{\ve{a}}])^H] \\
%			&= \m{T} \m{C}_{\ve{a}_{\mathbb{R}}\ve{a}_{\mathbb{R}}} \m{T}^H \label{equ:CWCU_Journal006a}\\
			&= \begin{bmatrix} \m{C}_{\ve{a}\ve{a}} & \tilde{\m{C}}_{\ve{a}\ve{a}} \\ 
				 \tilde{\m{C}}_{\ve{a}\ve{a}}^* & \m{C}_{\ve{a}\ve{a}}^* \end{bmatrix} = \underline{\m{C}}_{\ve{a}\ve{a}}^H \in\mathbb{C}^{2N_\ve{a}\times 2N_\ve{a}}, 
				 \label{equ:CWCU_Journal006}
\end{align}
with $\m{C}_{\ve{a}\ve{a}} = E_{\ve{a}}[(\ve{a} - E_{\ve{a}}[\ve{a}])(\ve{a} - E_{\ve{a}}[\ve{a}])^H]$ as the (Hermitian and positive semi-definite) covariance matrix and $\tilde{\m{C}}_{\ve{a}\ve{a}} = E_{\ve{a}}[(\ve{a} - E_{\ve{a}}[\ve{a}])(\ve{a} - E_{\ve{a}}[\ve{a}])^T]$ as the complementary covariance matrix. For $\m{C}_{\ve{a}\ve{a}}$ and $\tilde{\m{C}}_{\ve{a}\ve{a}}$ we have $\m{C}_{\ve{a}\ve{a}} =\m{C}_{\ve{a}\ve{a}}^H$ and $\tilde{\m{C}}_{\ve{a}\ve{a}} = \tilde{\m{C}}_{\ve{a}\ve{a}}^T$, respectively.

$\tilde{\m{C}}_{\ve{a}\ve{a}}$ is sometimes also referred to as pseudo-covariance matrix or conjugate covariance matrix. If $\tilde{\m{C}}_{\ve{a}\ve{a}} = \m{0}$, then the vector $\ve{a}$ is called \textit{proper}, otherwise 
\textit{improper} \cite{On_the_Circularity_of_a_Complex_Random_Variable, Essential_Statistics_and_Tools_for_Complex_Random_Variables, Second-order_analysis_of_improper_complex_random_vectors_and_processes, A_Complex_Generalized_Gaussian_Distribution_Characterization_Generation_and_Estimation}. For scalar Gaussian random variables, properness means that the equipotential lines of its PDF plotted in the complex plane are circles. If those equipotential lines are elliptical, then the scalar Gaussian random variable is improper. If $\ve{a}$ is real valued, then $\m{C}_{\ve{a}\ve{a}} = \tilde{\m{C}}_{\ve{a}\ve{a}}$. 
%General scalar random variables are called proper if the variances of the real and imaginary parts are equal. 

 Let $\ve{x}$ be the parameter vector to be estimated and $\ve{y}$ be the measurement vector, then a general widely linear estimator takes on the form in \eqref{equ:WB_Real002} (or an affine version of \eqref{equ:WB_Real002}).
%\begin{equation}
%	\hat{\ve{x}} = \m{E} \ve{y} + \m{F} \ve{y}^{*} + \ve{b}. \label{equ:prel_010}
%\end{equation}
In general, widely linear estimators are superior to their linear counterparts as soon as the measurements $\ve{y}$ turn improper. Applications for widely linear estimators are investigated in \cite{Statistical-Signal-Processing-of-Complex-Valued-Data-The-Theory-of-Improper-and-Noncircular-Signals, Widely-and-semi-widely-linear-processing-of-quaternion-vectors, Optimal-widely-linear-MVDR-beamforming-for-noncircular-signals, Component-Wise_Conditionally_Unbiased_Widely_Linear_MMSE_Estimation}. Another way to express the estimator in \eqref{equ:WB_Real002} is its augmented version
\begin{equation}
	\underline{\hat{\ve{x}}} = \begin{bmatrix}  \m{E} & \m{F} \\ \m{F}^{*} & \m{E}^{*} \end{bmatrix} 
														 \begin{bmatrix} \ve{y} \\ \ve{y}^{*} \end{bmatrix}
													 = \underline{\m{E}}\, \underline{\ve{y}} .\label{equ:WB_Real004d}
\end{equation}
This work focusses on the BLUE \cite{Kay-Est.Theory}, the BWLUE \cite{Statistical-Signal-Processing-of-Complex-Valued-Data-The-Theory-of-Improper-and-Noncircular-Signals} and the WLMMSE estimator \cite{Widely-linear-estimation-with-complex-data, Statistical-Signal-Processing-of-Complex-Valued-Data-The-Theory-of-Improper-and-Noncircular-Signals} for zero mean parameter vectors given in \eqref{equ:WB_Real004b}--\eqref{equ:CWCU_Journal008aa}, respectively. 
\begin{equation}
\hat{\ve{x}}_\text{B} = \left(\m{H}^H \m{C}_{\ve{n}\ve{n}}^{-1}\,\m{H} \right)^{-1} \m{H}^H \m{C}_{\ve{n}\ve{n}}^{-1}\,\ve{y}, \label{equ:WB_Real004b}
\end{equation}
\begin{equation}
\hat{\vea{x}}_\text{WB} = \left(\ma{H}^H \ma{C}_{\ve{n}\ve{n}}^{-1}\,\ma{H} \right)^{-1} \ma{H}^H \ma{C}_{\ve{n}\ve{n}}^{-1}\,\vea{y}, \label{equ:WB_Real004c}
\end{equation}
\begin{equation}
	\underline{\hat{\ve{x}}}_\text{WL} =  \underline{\m{C}}_{\ve{x}\ve{x}}\,\underline{\m{H}}^H \left(\underline{\m{H}}\, \underline{\m{C}}_{\ve{x}\ve{x}}\,\underline{\m{H}}^H + \underline{\m{C}}_{\ve{n}\ve{n}}\right)^{-1}\underline{\ve{y}}, \label{equ:CWCU_Journal008aa}
\end{equation} 
where
\begin{equation}
\ma{H} = \begin{bmatrix}
\m{H} & \m{0} \\
\m{0} & \m{H}^*
\end{bmatrix}.
\end{equation}
Note that the BWLUE reduces to the BLUE for proper noise.

\section{Derivation} \label{sec:Derivation}

In the following, the derivation of the BWLUE for real valued parameter vectors is given, assuming proper noise statistics ($\tilde{\m{C}}_{\ve{n}\ve{n}} = \m{0}$). An extension to improper noise will be handled in an upcoming paper. We start with some arguments about the BWLUE. The BWLUE for the $i^{\text{th}}$ element of the parameter vector written in the form
\begin{align}
\hat{x}_i =& \ve{e}_i^H\ve{y} + \ve{f}_i^H\ve{y}^* \label{equ:WB_Real005} \\
=& \begin{bmatrix} \ve{e}_i^H & \ve{f}_i^H  \end{bmatrix} \begin{bmatrix} \ve{y} \\ \ve{y}^*  \end{bmatrix}  \label{equ:WB_Real006} \\
=& \ve{w}_i^H \vea{y} \label{equ:WB_Real007}
\end{align}
can be derived by minimizing the cost function \cite{Statistical-Signal-Processing-of-Complex-Valued-Data-The-Theory-of-Improper-and-Noncircular-Signals}
\begin{equation}
J = \ve{w}_i^H \ma{C}_{\ve{n}\ve{n}} \ve{w}_i \label{equ:WB_Real008}
\end{equation}
s.t. the usual unbiasedness constraint $E_{\ve{n}}[\hat{x}_i]=x_i$ which simply follows to
\begin{equation}
\ve{w}_i^H \ma{H} = \ve{u}_i^T,\label{equ:WB_Real009}
\end{equation}
where $\ve{u}_i^T$ is a zero row vector of size $1 \times 2N_\ve{x}$ with a '1' at its $i^\text{th}$ position. 
To obtain our targeted BWLUE for real parameter vectors we in contrast to the ordinary BWLUE enforce
\begin{align}
\im{\hat{x}_i}& = 0  \label{equ:WB_Real009g} \\
E_{\ve{n}}[\re{\hat{x}_i}] = E_{\ve{n}}[\hat{x}_i]&=x_i, \label{equ:WB_Real009h}
\end{align}
where $\re{ \cdot }$ and $\im{ \cdot }$ denote the real and imaginary part, respectively. From \eqref{equ:WB_Real009g} one can easily show that the choice $\ve{e}_i^H = \ve{f}_i^T$ is necessary and sufficient to make $\hat{x}_i$ real valued independent of the concrete realization of $\ve{y}$. Incorporating this result into \eqref{equ:WB_Real009h} leads to
\begin{align}
E_{\ve{n}}[ \hat{x}_i ] =& E_{\ve{n}}\left[ \ve{e}_i^H \ve{y} + \ve{e}_i^T \ve{y}^*   \right] \label{equ:WB_Real011} \\
=&\ve{e}_i^H  \m{H} \ve{x} + \ve{e}_i^T \m{H}^* \ve{x} \label{equ:WB_Real013} \\
=& \left(\ve{e}_i^H  \m{H} + \ve{e}_i^T \m{H}^*\right) \ve{x}, \label{equ:WB_Real014} 
\end{align}
hence the unbiased constraint for this estimator is
\begin{equation}
\ve{e}_i^H  \m{H} + \ve{e}_i^T \m{H}^* = \ve{u}_i^T,\label{equ:WB_Real015}
\end{equation}
with $\ve{u}_i^T$ being a zero row vector of size $1 \times N_\ve{x}$ with a '1' at its $i^\text{th}$ position. Note that \eqref{equ:WB_Real015} is a less stringent requirement than \eqref{equ:WB_Real009}. Altogether this leads to the constrained optimization problem
\begin{align}
&\ve{e}_i = \text{argmin}\left( \begin{bmatrix} \ve{e}_i^H & \ve{e}_i^T  \end{bmatrix} \ma{C}_{\ve{n}\ve{n}} \begin{bmatrix} \ve{e}_i \\ \ve{e}_i^*  \end{bmatrix} \right) \nonumber \\
& \hspace{5mm} \text{s.t.} \hspace{5mm} \ve{e}_i^H \m{H} + \ve{e}_i^T \m{H}^* = \ve{u}_i^T, \label{equ:WB_Real010}
\end{align}
which can be solved by utilizing the Lagrange multiplier method. The Lagrange cost function follows to
\begin{align}
J' =&  \begin{bmatrix} \ve{e}_i^H & \ve{e}_i^T  \end{bmatrix} \ma{C}_{\ve{n}\ve{n}} \begin{bmatrix} \ve{e}_i \\ \ve{e}_i^*  \end{bmatrix}  + \lambda^T \left( \m{H}^H \ve{e}_i +  \m{H}^T \ve{e}_i^* - \ve{u}_i  \right) \nonumber \\
=&   \ve{e}_i^H \m{C}_{\ve{n}\ve{n}} \ve{e}_i  + \ve{e}_i^T  \m{C}_{\ve{n}\ve{n}}^*  \ve{e}_i^* + \lambda^T \left( \m{H}^H \ve{e}_i +  \m{H}^T \ve{e}_i^*  - \ve{u}_i \right). \nonumber
\end{align}
Note that the Lagrange multiplier $\lambda$ is real valued since the constraint is real valued. Taking the partial derivative of $J'$ w.r.t. $\ve{e}_i^*$ (using Wirtinger's calculus \cite{Wirtinger1927}, i.e. treating $\ve{e}_i$ as static) results in
\begin{equation}
\frac{\partial J'}{\partial \ve{e}_i^*} =  \m{C}_{\ve{n}\ve{n}} \ve{e}_i + \m{C}_{\ve{n}\ve{n}} \ve{e}_i + \m{H} \lambda.  \label{equ:WB_Real018}
\end{equation}
Setting this result equal to zero allows to identify
\begin{equation}
\ve{e}_i = - \frac{1}{2}\m{C}_{\ve{n}\ve{n}}^{-1}  \m{H} \lambda.  \label{equ:WB_Real019}
\end{equation}
Inserting \eqref{equ:WB_Real019} into the constraint in \eqref{equ:WB_Real010} yields
%\begin{equation}
%- \frac{1}{2}\lambda^T \m{H}^H\m{C}_{\ve{n}\ve{n}}^{-1} \m{H} - \frac{1}{2}\lambda^T \m{H}^T  \left(\m{C}_{\ve{n}\ve{n}}^{-1}\right)^* \m{H}^* = \ve{u}_i^T \label{equ:WB_Real020}
%\end{equation}
\begin{equation}
\lambda^T = -2\ve{u}_i^T \left(\m{H}^H\m{C}_{\ve{n}\ve{n}}^{-1} \m{H} +  \m{H}^T  \left(\m{C}_{\ve{n}\ve{n}}^{-1}\right)^* \m{H}^*  \right)^{-1}.  \label{equ:WB_Real021a} 
\end{equation}
Inserting \eqref{equ:WB_Real021a} into \eqref{equ:WB_Real019} leads to
%\begin{equation}
%\ve{e}_i = \m{C}_{\ve{n}\ve{n}}^{-1}  \m{H} \m{M}^{-1} \ve{u}_i  \label{equ:WB_Real025}
%\end{equation}
\begin{align}
\ve{e}_i^H =& \ve{u}_i^T \left(\m{H}^H\m{C}_{\ve{n}\ve{n}}^{-1} \m{H} +  \m{H}^T  \left(\m{C}_{\ve{n}\ve{n}}^{-1}\right)^* \m{H}^*  \right)^{-1}  \m{H}^H\m{C}_{\ve{n}\ve{n}}^{-1}.   \label{equ:WB_Real027}
\end{align}
Combining \eqref{equ:WB_Real027} with \eqref{equ:WB_Real003} results in
\begin{align}
&\hat{x}_i = \ve{u}_i^T \left(\m{H}^H\m{C}_{\ve{n}\ve{n}}^{-1} \m{H} +  \m{H}^T  \left(\m{C}_{\ve{n}\ve{n}}^{-1}\right)^* \m{H}^*  \right)^{-1}  \m{H}^H\m{C}_{\ve{n}\ve{n}}^{-1}\ve{y} \nonumber \\
& + \ve{u}_i^T \left(\m{H}^H\m{C}_{\ve{n}\ve{n}}^{-1} \m{H} +  \m{H}^T  \left(\m{C}_{\ve{n}\ve{n}}^{-1}\right)^* \m{H}^*  \right)^{-1}  \m{H}^T\left(\m{C}_{\ve{n}\ve{n}}^{-1}\right)^*\ve{y}^*  \label{equ:WB_Real029}
\end{align}
The vector estimator immediately follows to
\begin{equation}
\hat{\ve{x}} = \m{E}\ve{y} + \m{E}^*\ve{y}^*, \label{equ:WB_Real061a} 
\end{equation}
where the estimator matrix $\m{E}$ follows from \eqref{equ:WB_Real027}  as
\begin{align}
\m{E} =& \left(\m{H}^H\m{C}_{\ve{n}\ve{n}}^{-1} \m{H} +  \m{H}^T  \left(\m{C}_{\ve{n}\ve{n}}^{-1}\right)^* \m{H}^*  \right)^{-1}  \m{H}^H\m{C}_{\ve{n}\ve{n}}^{-1} \label{equ:WB_Real061b} 
%\m{E}^* =&\left(\m{H}^H\m{C}_{\ve{n}\ve{n}}^{-1} \m{H} +  \m{H}^T  %\left(\m{C}_{\ve{n}\ve{n}}^{-1}\right)^* \m{H}^*  \right)^{-1}  \m{H}^T%\left(\m{C}_{\ve{n}\ve{n}}^{-1}\right)^*. \label{equ:WB_Real061c}
\end{align}
This is the final result for the BWLUE for real valued parameter vectors and proper noise.

\section{Discussion} \label{sec:Discussion}

In a first step we verify that this result yields an unbiased estimator. The expectation of $\hat{\ve{x}}$ follows to
\begin{align}
&E_{\ve{n}}[\hat{\ve{x}}] =  \left(\m{H}^H\m{C}_{\ve{n}\ve{n}}^{-1} \m{H} +  \m{H}^T  \left(\m{C}_{\ve{n}\ve{n}}^{-1}\right)^* \m{H}^*  \right)^{-1}  \m{H}^H\m{C}_{\ve{n}\ve{n}}^{-1}\m{H}\ve{x} \nonumber \\
&  +  \left(\m{H}^H\m{C}_{\ve{n}\ve{n}}^{-1} \m{H} +  \m{H}^T  \left(\m{C}_{\ve{n}\ve{n}}^{-1}\right)^* \m{H}^*  \right)^{-1}  \m{H}^T\left(\m{C}_{\ve{n}\ve{n}}^{-1}\right)^*\m{H}^*\ve{x} \nonumber \\
%&= \left(\m{H}^H\m{C}_{\ve{n}\ve{n}}^{-1} \m{H} +  \m{H}^T  \left(\m{C}_{\ve{n}\ve{n}}^{-1}\right)^* \m{H}^*  \right)^{-1} \left(\m{H}^H\m{C}_{\ve{n}\ve{n}}^{-1} \m{H} +  \m{H}^T  \left(\m{C}_{\ve{n}\ve{n}}^{-1}\right)^* \m{H}^*  \right)\ve{x} \label{equ:WB_Real033} \nonumber \\
& \hspace{6mm} = \ve{x}. 
\end{align}
which proofs that this estimator is unbiased. 

Another formulation of \eqref{equ:WB_Real061a}-\eqref{equ:WB_Real061b} can be derived as
\begin{align}
&\hat{\ve{x}} =  
%\left(\m{H}^H\m{C}_{\ve{n}\ve{n}}^{-1} \m{H} +  \m{H}^T  \left(\m{C}_{\ve{n}\ve{n}}^{-1}\right)^* \m{H}^*  \right)^{-1}  \m{H}^H\m{C}_{\ve{n}\ve{n}}^{-1}\ve{y} %\nonumber \\
%& \hspace{2mm} +  \left(\m{H}^H\m{C}_{\ve{n}\ve{n}}^{-1} \m{H} +  \m{H}^T  \left(\m{C}_{\ve{n}\ve{n}}^{-1}\right)^* \m{H}^*  \right)^{-1}  \m{H}^T\left(\m{C}%_{\ve{n}\ve{n}}^{-1}\right)^*\ve{y}^* \nonumber \\
% & \hspace{2mm} =  
\left( 2\re{ \m{H}^H\m{C}_{\ve{n}\ve{n}}^{-1} \m{H}}  \right)^{-1} \left( \m{H}^H\m{C}_{\ve{n}\ve{n}}^{-1}\ve{y} +  \m{H}^T\left(\m{C}_{\ve{n}\ve{n}}^{-1}\right)^*\ve{y}^* \right) \nonumber \\
  & \hspace{2mm} =  \left( \re{ \m{H}^H\m{C}_{\ve{n}\ve{n}}^{-1} \m{H}}  \right)^{-1}\re{ \m{H}^H\m{C}_{\ve{n}\ve{n}}^{-1}\ve{y} }. \label{equ:WB_Real036} 
\end{align}
This is the most compact formulation of the BWLUE for real valued parameter vectors. 

Assuming the special case where the term $\m{H}^H\m{C}_{\ve{n}\ve{n}}^{-1} \m{H}$ is real valued we obtain
\begin{align}
\hat{\ve{x}} =&   \left( \re{ \m{H}^H\m{C}_{\ve{n}\ve{n}}^{-1} \m{H}}  \right)^{-1}\re{ \m{H}^H\m{C}_{\ve{n}\ve{n}}^{-1}\ve{y} } \label{equ:WB_Real037}  \\
=& \left(  \m{H}^H\m{C}_{\ve{n}\ve{n}}^{-1} \m{H}  \right)^{-1}\re{ \m{H}^H\m{C}_{\ve{n}\ve{n}}^{-1}\ve{y} } \label{equ:WB_Real038} \\
=& \re{ \left(  \m{H}^H\m{C}_{\ve{n}\ve{n}}^{-1} \m{H}  \right)^{-1}  \m{H}^H\m{C}_{\ve{n}\ve{n}}^{-1}\ve{y}  }. \label{equ:WB_Real039}
\end{align}
In that case, the BWLUE for real parameter vectors coincides with the real part of the BLUE in \eqref{equ:WB_Real004b}. Furthermore, it also coincides with the real part of the BWLUE in \eqref{equ:WB_Real004c} since the noise is assumed to be proper. 
%The BWLUE for real valued parameter vectors can be brought in augmented form according to \eqref{equ:WB_Real004d} where $\m{F}$ has to be replaced by $\m{E}^*$. Then, the augmented error covariance matrix can be derived as
%\begin{equation}
%\ma{C}_{\ve{e}\ve{e}} = \ma{E} \ma{C}_{\ve{n}\ve{n}} \ma{E}^H. \label{equ:WB_Real066}
%\end{equation}
%The ordinary error covariance matrix $\m{C}_{\ve{e}\ve{e}}$ is the NW-block of $\ma{C}_{\ve{e}\ve{e}}$. 

The covariance matrix of $\hat{\ve{x}}$ can be shown to be
\begin{align}
\m{C}_{\hat{\ve{x}}\hat{\ve{x}}} =& \m{E} \,\m{C}_{\ve{n}\ve{n}}\, \m{E}^H + \m{E}^* \m{C}_{\ve{n}\ve{n}}^* \,\m{E}^T \label{equ:WB_Real067} \\
=& 2\, \re{\m{E} \m{C}_{\ve{n}\ve{n}} \m{E}^H}  \label{equ:WB_Real067_1} \\
=& \left(\m{H}^H\m{C}_{\ve{n}\ve{n}}^{-1} \m{H} +  \m{H}^T  \left(\m{C}_{\ve{n}\ve{n}}^{-1}\right)^* \m{H}^*  \right)^{-1}.\label{equ:WB_Real067_2}
\end{align}
The expression for the BWLUE for real valued parameter vectors in \eqref{equ:WB_Real061a}-\eqref{equ:WB_Real061b} could also have been derived by minimizing the Bayesian mean square error (BMSE) cost function $E_{\ve{y},\ve{x}}[|\hat{x}_i - x_i|^2]$ subject to the same constraint as in \eqref{equ:WB_Real010}, such that this estimator can also be interpreted in a Bayesian sense. In this sense, the Bayesian error covariance matrix $\m{C}_{\ve{e}\ve{e}}$, where $\ve{e} = \hat{\ve{x}} - \ve{x}$, is equal to $\m{C}_{\hat{\ve{x}}\hat{\ve{x}}}$ in \eqref{equ:WB_Real067}--\eqref{equ:WB_Real067_2} and the minimum BMSEs can be found on the main diagonal of $\m{C}_{\ve{e}\ve{e}}$.

For every complex linear or widely linear estimator such as the WLMMSE estimator, there exists an equivalent real valued estimator derived from an extended real valued linear model. The formulation as complex linear model and according complex estimator has the advantage of being much more compact. The equivalent real valued model to \eqref{equ:WB_Real001} is
\begin{equation}
\begin{bmatrix}
\re{\ve{y}} \\ \im{\ve{y}}
\end{bmatrix} = \begin{bmatrix}
\re{\m{H} } \\ \im{\m{H} }
\end{bmatrix}  \ve{x} + \begin{bmatrix}
\re{\ve{n}} \\ \im{\ve{n}}
\end{bmatrix}. \label{equ:WB_Real001a}
\end{equation}
The BLUE applied on this real valued linear model follows to
\begin{align}
&\hat{\ve{x}} = \nonumber \\
&\left(\begin{bmatrix}
\re{\m{H} } \\ \im{\m{H} }
\end{bmatrix}^H \bar{\m{C}}_{\ve{n}\ve{n}}^{-1}\begin{bmatrix}
\re{\m{H} } \\ \im{\m{H} }
\end{bmatrix} \right)^{-1} \begin{bmatrix}
\re{\m{H} } \\ \im{\m{H} }
\end{bmatrix}^H \bar{\m{C}}_{\ve{n}\ve{n}}^{-1}\begin{bmatrix}
\re{\ve{y}} \\ \im{\ve{y}}
\end{bmatrix} \label{equ:WB_Real004e}
\end{align}
where $\bar{\m{C}}_{\ve{n}\ve{n}}$ is an appropriately modified covariance matrix of the noise vector in \eqref{equ:WB_Real001a}. 
Eq.~\eqref{equ:WB_Real004e} is equivalent to \eqref{equ:WB_Real061a}-\eqref{equ:WB_Real061b} and \eqref{equ:WB_Real036}. However, comparing the matrices to be inverted in \eqref{equ:WB_Real036} and in \eqref{equ:WB_Real004e} reveals the advantage of the proposed estimator.  Hence, applying the proposed estimator allows for a significant complexity reduction.

\section{Simulation Results} \label{sec:Simulation Results}

To compare the performance of the derived estimator we consider the case where a discrete-time approximation of a real valued analogue impulse response is estimated based on spectral measurements of an analogue linear and time invariant (LTI) system. The measurements are performed at equidistant frequency steps $f = k\Delta f$ with $k=0,\cdots, 19$. The complex measurements assembled in a vector are denoted as $\ve{y}$. These measurements represent the single sided frequency response of the analogue LTI system. An extension to a double sided frequency response causes the measurement noises to be perfectly correlated. This would further result in a singular $\m{C}_{\ve{n}\ve{n}}$, preventing the application of the BLUE. The discrete-time impulse response $\ve{x}$, which is assumed to have a length of $N_\ve{x} = 5$, is connected with the measurements via
\begin{equation}
\ve{y} = T_S \m{H}\ve{x} + \ve{n}. \label{equ:WB_Real070}
\end{equation}
In \eqref{equ:WB_Real070}, $\m{H}$ is given by the first 5 columns and the first 20 rows of a DFT matrix of size $40 \times 40$, and $T_S$ is the sampling time of $\ve{x}$ and is necessary in order to combine analogue spectral measurements with a discrete-time impulse response \cite{Discrete-time-Signal-Processing}. In the following, $T_S$ is set to 1 for simplicity.

In order to compare the BWLUE for real valued parameter vectors with the Bayesian WLMMSE estimator, the simulations are carried out in a Bayesian framework where the impulse response $\ve{x}$ is randomly generated with $E[\ve{x}] = \ve{0}$ and $\m{C}_{\ve{x}\ve{x}} = \tilde{\m{C}}_{\ve{x}\ve{x}} = \m{I}$. 
$\ve{n}$ in \eqref{equ:WB_Real070} is assumed to be a zero mean and proper noise vector with covariance matrix $\m{C}_{\ve{n}\ve{n}} = \sigma_\ve{n}^2 \m{I}$, and $\sigma_\ve{n}^2$ is varied between $10^{-3}$ and $10^{2}$. 

\begin{figure}[tb]
\begin{center}
\begin{tikzpicture}
\begin{loglogaxis}[compat=newest, 
width=0.9\columnwidth, height = .6\columnwidth, xlabel=$\sigma_\ve{n}^2$, 
ylabel style={align=center}, 
ylabel style={text width=3.4cm},
ylabel={average BMSE}, 
legend pos=north east, 
legend cell align=left,
legend columns=2, 
        legend style={
                    % the /tikz/ prefix is necessary here...
                    % otherwise, it might end-up with `/pgfplots/column 2`
                    % which is not what we want. compare pgfmanual.pdf
            /tikz/column 2/.style={
                column sep=5pt,
            },
        font=\small},
xmin = 0.001,
xmax = 100,
ymax = 500,
ymin = 0.00001,
grid=major,
%restrict y to domain=-inf:20,
legend style={
at={(-0.25,1.3)},
anchor=north west}
]

% BLUE
\addplot[line width=1pt, color=blue] table[x index =1, y index =3] {BWLUE_Real_ICASSP.dat};
\addlegendentry{{BLUE}}

% Real part of BLUE
\addplot[line width=1pt, color=blue, style=dashed] table[x index =1, y index =5] {BWLUE_Real_ICASSP.dat};
\addlegendentry{{Real part of the BLUE}}

% WLMMSE estimato
\addplot[line width=1pt, color=red, style=solid] table[x index =1, y index =2] {BWLUE_Real_ICASSP.dat};
\addlegendentry{{WLMMSE est.}}

% BWLUE for real parameter vectors
\addplot[line width=1pt, color=green] table[x index =1, y index =4] {BWLUE_Real_ICASSP.dat};
\addlegendentry{{BWLUE for real parameter vectors}}

\end{loglogaxis}
\end{tikzpicture}
\caption{Average BMSEs of the estimated impulse response for the BLUE in \eqref{equ:WB_Real004b}, the real part of the BLUE according to \eqref{equ:WB_Real039}, the WLMMSE estimator in \eqref{equ:CWCU_Journal008aa} and the BWLUE for real valued parameter vectors in \eqref{equ:WB_Real061a}-\eqref{equ:WB_Real061b}. The BLUE and the BWLUE in \eqref{equ:WB_Real004c} coincide since the noise $\ve{n}$ is proper. \label{fig:Example1} }
\end{center}
\end{figure}
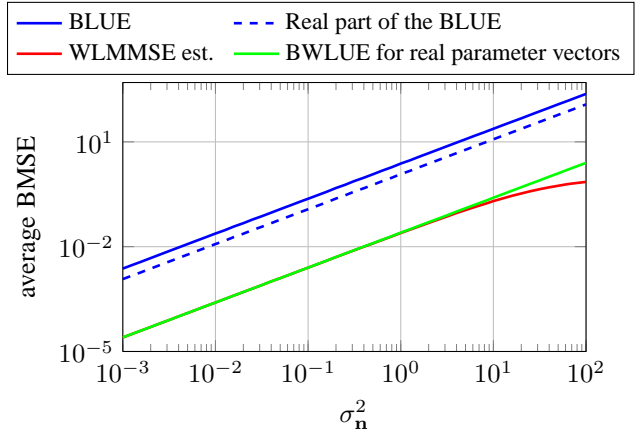

The resulting average BMSE values, defined as the average BMSEs of the vector estimator's elements, are presented in Fig.~\ref{fig:Example1}. There, the performance of the BLUE as well as the real part of the BLUE as in \eqref{equ:WB_Real039} is shown. By taking only the real part of the BLUE, the performance can be slightly increased. The BLUE and the BWLUE coincide since the noise is proper, this is why only the results of the BLUE are shown. A significant increase in performance is provided by the WLMMSE estimator, but the derived BWLUE for real parameter vectors which does not use prior statistical knowledge about the parameters almost reaches the performance of the WLMMSE estimator for noise variances up to $10^1$. In this area, the prior knowledge in form of $\m{C}_{\ve{x}\ve{x}}$ and $\tilde{\m{C}}_{\ve{x}\ve{x}}$ is weak compared to the information content of the measurements. For larger noise variances, the prior knowledge begins to increase the performance of the WLMMSE estimator compared to the BWLUE for real valued parameter vectors. However, for all values of $\sigma_\ve{n}^2$, the proposed estimator outperforms the classical BLUE and BWLUE by approximately two orders of magnitude in BMSE.

\section{Conclusion} \label{sec:Conclusion}

In this work the task of estimating a real valued parameter vector based on complex measurements has been investigated. This task has been analyzed in the classical sense, i.e. for the assumption of a deterministic but unknown parameter vector. A novel classical estimator has been derived which can incorporate the knowledge that the parameter vector is real valued. This estimator outperforms other classical estimators such as the BLUE or the BWLUE. A simulation example was shown revealing that the proposed estimator can compete with the Bayesian WLMMSE estimator, which requires first and second order statistics about the parameter vector. Proper noise statistics have been assumed for the derivation. An extension to improper noise will be handled in upcoming work.

%% Below is an example of how to insert images. Delete the ``\vspace'' line,
%% uncomment the preceding line ``\centerline...'' and replace ``imageX.ps''
%% with a suitable PostScript file name.
%% -------------------------------------------------------------------------
%\begin{figure}[htb]
%
%\begin{minipage}[b]{1.0\linewidth}
%  \centering
%  \centerline{\includegraphics[width=8.5cm]{image1}}
%%  \vspace{2.0cm}
%  \centerline{(a) Result 1}\medskip
%\end{minipage}
%%
%\begin{minipage}[b]{.48\linewidth}
%  \centering
%  \centerline{\includegraphics[width=4.0cm]{image3}}
%%  \vspace{1.5cm}
%  \centerline{(b) Results 3}\medskip
%\end{minipage}
%\hfill
%\begin{minipage}[b]{0.48\linewidth}
%  \centering
%  \centerline{\includegraphics[width=4.0cm]{image4}}
%%  \vspace{1.5cm}
%  \centerline{(c) Result 4}\medskip
%\end{minipage}
%%
%\caption{Example of placing a figure with experimental results.}
%\label{fig:res}
%%
%\end{figure}

%\vfill
%\pagebreak


\bibliography{References}

%\begin{thebibliography}{20}
%
%\bibitem{Schreier-Buch} P.J. Schreier, L.L. Scharf; ''Statistical signal processing of complex-valued data: the theory of improper and noncircular signals,'' 2010, Cambridge University Press.
%
%\bibitem{Schreier-2011} T. Adali, P.J. Schreier, L.L. Scharf; ''Complex-Valued Signal Processing: The Proper Way to Deal With Impropriety,'' In \emph{IEEE Trans. Signal Process.}, vol. 59, issue 11, pp. 5101-5125, 2011.
%%\newline doi: 10.1109/TSP.2011.2162954
%
%
%	
%\end{thebibliography}

\end{document}